\documentclass[twoside, leqno]{article}

\usepackage{fullpage}

\usepackage{ltexpprt}
\usepackage{hyperref}

\usepackage{graphicx}
\usepackage{amsmath,amsfonts,amssymb}
\usepackage{fullpage}
\usepackage{url}
\usepackage{todonotes}

\newtheorem{definition}{Definition}
\newcommand{\qed}{\Box}

\usepackage{xcolor}

\begin{document}

\title{\Large On the size of the neighborhoods of a word}
\author{Cedric Chauve\thanks{Department of Mathematics, Simon Fraser University, Canada. cedric.chauve@sfu.ca}
\and Louxin Zhang\thanks{Department of Mathematics, National University of Singapore. matzlx@nus.edu.sg}}

\date{}

\maketitle

\begin{abstract} \small\baselineskip=9pt 
The $d$-neighborhood of a word $w$ in the Levenshtein distance is the set of all words at distance at most $d$ from $w$. 
Generating the neighborhood of a word $w$, or related sets of words such as the condensed neighborhood or the super-condensed neighborhood has applications in the design of approximate pattern matching algorithms. 
It follows that bounds on the maximum size of the neighborhood for the words of a given length can be used in the complexity analysis of such approximate pattern matching algorithms.
In this note, we present exact formulas for the sizes of the condensed and super condensed neighborhoods of unary words, {establish a novel upper bound  and  prove a conjectured upper bound for the size of the condensed neighborhoods of an arbitrary word.}

\end{abstract}

%% ----------------------------------------------------------------------------------------------
\section{Introduction}
\label{sec:introduction}

%Paragraph 1 [Motivation]: sequence search in bioinformatics.
Aiming to  search for all approximate occurrences of a query sequence within a text,
a problem known as approximate pattern matching, is at the heart of many basic applications in bioinformatics~\cite{DBLP:books/cu/Gusfield1997,DBLP:journals/csur/Navarro01}, particularly in searching large biological sequence databases with the BLAST tool~\cite{BLAST}. 
The BLAST algorithm proceeds in two phases. 
First, it identifies \emph{seeds} that are short sequences present in both the query sequence and the target text.  These seed occurrences in the searched { text} are then \emph{extended}, using dynamic programming, to form approximate occurrences of the query.
This approach to approximate pattern matching is known as the \emph{seed-and-extend} approach.

%Paragraph 2 [Main concept]: neighborhoods and their use in sequence search.
Conceptually, a key part of the first phase of seed-and-extend methods { is} generating, for every subsequence $w$ of length $n$ of the query and for a distance value $d$ 
%-- chosen { according} to the maximum distance between the query and the searched approximate occurrences of the query -- 
the set of all words at Levenshtein distance (also known as edit distance) from $w$ at most $d$, known as the \emph{$d$-neighborhood} of $w$.
{ The words of the neighborhood are then used as seeds.}
In practice, variants of the neighborhood concept are used in approximate pattern matching algorithms{ .
For example, BLAST uses the condensed neighborhood that excludes words having a prefix that is itself a word in the neighborhood~\cite{DBLP:books/daglib/p/Myers13}, motivated by the property that any prefix of a seed is a seed itself. }
%An example is the condensed neighborhood, used in BLAST~\cite{DBLP:books/daglib/p/Myers13}, which is obtained by discarding from the neighborhood any word having a prefix that is itself a word in the neighborhood. 
Another variant is the super condensed neighborhood, which discards words having a subword already in the neighborhood~\cite{DBLP:journals/jda/RussoO07}.

Bounds on the maximum size of the neighborhood over all words of a given length over a given alphabet play an important role in the complexity analysis of seed-and-extend approximate pattern matching algorithms. 
However, there are still few known results on this topic~\cite{DBLP:conf/cpm/Charalampopoulos20b, DBLP:conf/stringology/ChauveMP21, DBLP:journals/algorithmica/Myers94, DBLP:books/daglib/p/Myers13,DBLP:conf/lata/Touzet16}. 
{ For the size of the condensed neighborhood, motivated by the analysis of the complexity of BLAST, Myers provides in~\cite{DBLP:books/daglib/p/Myers13} a set of recurrences defining an upper-bound, and derives analytically from this recursions an upper-bound formula.
In~\cite{DBLP:conf/stringology/ChauveMP21}, the authors provide an asymptotic expression for the recurrences described in~\cite{DBLP:books/daglib/p/Myers13} and conjecture that the size of the condensed $d$-neighborhood of any word of length { {$n$}} over an alphabet of size $s$ is bounded above by $\frac{(2s-1)^d{ n}^d}{d!}$.}

%Paragraph 3 [Plan]: results and plan of the paper.
In this note, we provide several results on the size of the condensed and super condensed neighborhoods of a word.
In Section~\ref{sec:unary} we provide formulas for the size of the condensed and super condensed neighborhoods of unary words.
In Section~\ref{sec:arbitrary} we provide a { novel upper bound for} the size of the condensed neighborhood of arbitrary words 
{ and we use this formula} in Section~\ref{sec:conjecture} to prove the conjecture of~\cite{DBLP:conf/stringology/ChauveMP21}.
%% ----------------------------------------------------------------------------------------------
\section{Preliminaries}
\label{sec:preliminaries}

In this section, we introduce formal definitions and notations that will be used in this Note.

\smallskip\paragraph{Words}
Let $\Sigma$ be a finite set of characters, called an alphabet. 
{ A word $w$ over $\Sigma$ is an ordered sequence of characters $w=w_1w_2\cdots w_k$, where $w_i\in \Sigma$. 
Its length is defined as the number $k$ of characters appearing in $w$, denoted by $|w|$.
}
The empty word of length $0$ is denoted by $\epsilon$.
We denote by $\Sigma^{+}$ the set of all nonempty words and $\Sigma^{*}=\Sigma^{+}\cup \{\epsilon\}$.

For two words $u$ and $v$, we use $uv$ to denote the word obtained by concatenating $v$ and $v$. 
%Clearly, $u\epsilon=\epsilon u=u$.
We also define $u^0=\epsilon$ and, for a positive integer $k\geq 1$, $u^k$ is the concatenation of $k$ copies of $u$.
For two sets of words ${\cal  U}$ and ${\cal  V}$, we define ${\cal   UV} = \{uv \mid u\in {\cal U}, v\in {\cal V}\}$ as the set of concatenations of a word from ${\cal U}$ and a word from ${\cal V}$. 
A word $w$ is \textit{unary} if it consists of multiple occurrences of a single character from $\Sigma$, i.e., $w=\sigma^{|w|}$ for some $\sigma\in\Sigma$.

{ For words $u$ and $v$,  $u$ is said to be a prefix (resp. suffix) of $v$ if $v=uw$ (resp. $v=wu$) for some $w\in \Sigma^{*}$;  $u$ is said to be a subword of $v$ if $v=xuy$ for some $x,y\in \Sigma^{*}$.
}

\smallskip\paragraph{Sequence Alignment and Levenshtein distance}
The Levenshtein distance between two words $u$ and $v$, denoted by $d_{lev}(u,v)$, is the minimum number of edit operations that { are required} to transform $u$ into $v$, where edit operations can be:
\begin{itemize}
 \item Insertion: inserting a character at some position in a word;
 \item Deletion: deleting a character at some position from a word;
 \item Substitution:  replacing a character in a word with a different character.
\end{itemize}

{ An alignment $A$ between two words $u$ and $v$ on $\Sigma$ is a two-row array, where each row is a word on the alphabet $\Sigma \cup \{-\}$  and no column contains two occurrences of '$-$' such that the words obtained from the two rows by deleting all occurrences of '$-$' are $u$ (first row) and $v$ (second row), respectively. 
The character '$-$' is called a \textit{gap}. 
}
%We denote by $A_1$ (resp. $A_2$) the word on $\Sigma$ obtained by deleting gaps from the first (resp. second) row of $A$.
The {\it cost} of an alignment is the number of columns containing two different characters.
A column containing twice the same character is called a \textit{match} column.
A column with two different characters is a \textit{deletion} column { if} the bottom character is a gap, an \textit{insertion} column if the top character is a gap, and a \textit{mismatch} column otherwise.

The Levenshtein distance between two words $u$ and $v$ is equal to the minimum cost of an alignment $A$ between $u$ and $v$.
For example:
{\tt
\begin{center}
\begin{tabular}{llllll}
 a & l & - & i & g & n\\
 a & s & s & i & g & n
\end{tabular}
\end{center}
}
\noindent
represents an optimal (minimum cost) alignment between the words ``align" and ``assign". { The alignment  has cost $2$, as it contains  a substitution ($l$ to $s$) and  an insertion (an extra $s$ in the second row).
}
An optimal alignment between two words $u,v$, { therefore} their Levenshtein distance, can be computed in quadratic $O(|u||v|)$ time and space using dynamic programming~\cite{DBLP:books/cu/Gusfield1997}.

\smallskip\paragraph{Neighborhoods of a word}
For a non-negative integer $d$, the $d$-neighborhood of a query word $w$,  is the following subset of words:
\begin{equation}
\label{nbhood}
N(w, d)=\{x \in \Sigma^{\star} \mid d_{\rm lev}(x, w)\leq d\}.
\end{equation}
Clearly, each word of $N(w, d)$ contains at most $d+|w|$ characters and thus $N(w, d)$ is finite.

The {\it condensed $d$-neighborhood} of $w$, written as ${CN}(w, d)$, consists of the words of $N(w, d)$ which do not have a prefix in $N(w, d)$, that is, 
\begin{equation}
\label{c-nbhood}
CN(w, d)=N(w, d)\setminus \left[ N(w, d)\Sigma^{+}\right]. 
\end{equation}

Lastly, the {\it super-condensed  $d$-neighborhood} of $w$, written as $SCN(w, d)$, consists of the words of $N(w, d)$ that do not have a subword  in $N(w, d)$, that is, 
\begin{eqnarray}
\label{sc-nbhood}
 & &SCN(w, d)\\
 & =& N(w, d)\setminus \left[ \Sigma^{*}N(w, d)\Sigma^{+}
\cup \Sigma^{+}N(w, d)\Sigma^{*}\right].  \nonumber
\end{eqnarray}

%% ----------------------------------------------------------------------------------------------
\section{(Super)-Condensed Neighborhood for Unary Words}
\label{sec:unary}

In this section, we provide exact formulas for the size of the condensed and { super-condensed} neighborhoods of unary words.

\begin{proposition}
    Let $\Sigma$ be an alphabet consisting of $s$ characters and let $\sigma \in \Sigma$. 
    For any { positive integers $n$} and  $d$ such that 
    { $0< d< n$}, 
    \begin{equation}
    \label{eqn4}
        \vert CN({ \sigma^n}, d)\vert 
        = \sum_{n-d\leq m\leq n}
        {m-1\choose d+m-n}(s-1)^{d+m-n}.
    \end{equation}
    In particular, if $s=2$,
    \begin{equation}
    \label{eqn4_binary}
        \vert CN({ \sigma^n}, d)\vert={n\choose d}.
    \end{equation}
\end{proposition}
%\todo[inline]{k should be n in the subscript in (4). $w>=d>0$ to $0<d<w$, as reviewer pointed out d=w the combination number in the sum is not well-defined. Here, we changed w to n. A complete new proof as suggested by Reviewer 1}
\noindent {\bf Proof.} 
{ Assume that $m$ is a non-negative integer.
Let { $w=\sigma^n$,} $x=x_1x_2\cdots x_{m-1}x_m$ be a word on $\Sigma$ such that $x\in CN(w, d)$ and the alignment $A$:
\begin{center}
\begin{tabular}{ccccc}
$a_1$ & $a_2$ & $\cdots$ & $a_{s-1}$ & $a_s$ \\
$b_1$ & $b_2$ & $\cdots$ & $b_{s-1}$ & $b_s$ 
\end{tabular}
\end{center}
be an optimal alignment between $x$ and $w$, where $s\geq \max(m, n)$ and $x$ and $w$ appear in the first and second row, respectively.

First, $d_{\rm lev}(x, w)=d$. Otherwise, by the triangle inequality, 
$d_{\rm lev}(x_1x_2\cdots x_{m-1}, w)\leq d_{\rm lev}(x_1x_2\cdots x_{m-1}, x)+d_{\rm lev}(x, w)\leq 1+(d-1)\leq d$ and thus 
$x_1x_2\cdots x_{m-1}\in N(w, d)$, a contradiction with $x\in CN(w, d)$.

Second, $A$ does not contain any deletion and $x_m=\sigma$.
Assume $A$ contains at least one deletion and let the last deletion be the column ${a_j \brack b_j}$, that is $a_j=x_i$ for some $i$ and $b_j=-$ and $b_k=\sigma$ for $k=j+1, \cdots, s$. This implies that the alignment
\begin{center}
\begin{tabular}{cccccccc}
$a_1$ & $a_2$ & $\cdots$ &$a_j$ &$a_{j+1}$ & $\cdots$ & $a_{s-1}$ & $a_s$ \\
$b_1$ & $b_2$ & $\cdots$ & $b_{j+1}$ &$b_{j+2}$& $\cdots$ & $b_{s}$ & $-$ 
\end{tabular}
\end{center}
is also an optimal alignment between $x$ and $w$ and the first 
$s-1$ columns form also an optimal alignment between $x_1x_2\cdots x_{m-1}$ and $w$. 
Therefore,
$d_{\rm lev}(x_1x_2\cdots x_{m-1}, w)\leq d_{\rm lex}(x, w)-1\leq d-1$ and $x_1x_2\cdots x_{m-1}\in N(w, d)$, a contradiction. 
So the alignment $A$ does not contain any deletion column, and $m\leq n { =s}$.
If $x_m\neq\sigma$, then $x_1x_2\cdots x_{m-1}\in N(w, d)$, again a contradiction with $x \in CN(w, d)$.

In addition, since $m\leq n$ { and there is no deletion, $A$ contains  exactly} $|w|-|x|=n-m$ insertions. 
So $n-m \leq d_{\rm lev}(x, w) { =}  d$, implying that $m\geq n-d$.
}

%Let $A$ be an optimal alignment between a word $x=x_1\dots x_m$ of length $m$ and a unary word $w=\sigma^n$. Without loss of generality, {  we may assume that $x$ and $w$ appear on the row one and   row two of $A$, respectively.}
%In $A$, insertion, deletion and mismatch columns each incur a cost of $1$.
%Therefore, if $m\leq n$, $A$ does not contain any { deletion column ${x_i \choose -}$;} furthermore, insertion columns can only appear after the column containing $x_m$.
%{ This is for the following reasons.}

%Assume that $x\in CN(\sigma^w, d)$ and $m=|x|$. 
%If $m>w$,  for each $i>w$, $x_i$ must appear as ${x_i \choose -}$ in the alignment $A$. 
%Therefore, we have $d_{lev}(x_1x_2\cdots x_w, \sigma^w)\leq d$ and thus  $x_1x_2\cdots x_w\in N(\sigma^w, d)$, leading to a contradiction.We have thus proved that $m\leq k$.  Similarly, $x_m$ must appear in a match column. 

%Taken into account, both facts imply that the optimal alignment  between $x$ and $\sigma^w$ has the following structure:
%\begin{center}
%\begin{tabular}{ccccccc}
%$x_1$ & $x_2$ & $\cdots$ & $x_{m}$ & - & $\cdots $ & -\\
%$\sigma$ & $\sigma$ & $\cdots$ & $\sigma$ & $\sigma$ & $\cdots $ & $\sigma$ 
%\end{tabular}
%\end{center}

Since $d_{lev}(x, w)=d$ and $x_m=\sigma$, there are $d-(n-m)$ mismatch columns in the first { $n-1$ columns of $A$, so $x$ has exactly $d-(n-m)$ characters different from $\sigma$ in its first $m-1$ characters, the $n-d$ other characters of $x$ being $\sigma$. 
There are exactly ${m-1 \choose d-(n-m)} (s-1)^{d-(n-m)}$ such words. 
Consider two words $x$ and $y$ having the structure described above, with $x=x_1x_2\cdots x_{m}$, $y=y_1y_2\cdots y_{p}$, $m < p$, both $x$ and $y$ contain exactly $n-d$ occurrences of $\sigma$, and $x_m=y_p=\sigma$. 
}
Since $y_p=\sigma$, $y_1y_2\cdots y_m$ contains at most $n-d-1$ occurrences of $\sigma$, and $x$ cannot be a prefix of $y$.  

Therefore, in total, $CN(w, d)$ contains 
$$
\sum_{n-d\leq m\leq n}{m-1\choose d+m-n}(s-1)^{m+d-n}
$$ 
words, which proves Eqn. (\ref{eqn4}).

Substituting $s$ with $2$ in Eqn. (\ref{eqn4}),  we have
\begin{eqnarray*}
&& \vert CN(w, d)\vert \\
&=& {n-d-1 \choose 0} +  {n-d\choose 1} 
+ \cdots + {n-1\choose d}\\
&=& {n-d-1 +d+1\choose d} \;\;\;\;\;  \mbox{ (by the hockey stick identity)}\\
&=& {n\choose d}, 
\end{eqnarray*}
which proves Eqn. (\ref{eqn4_binary}). 
$\qed$
 
 % \begin{corollary}
 %  On the binary alphabet, for any $w$ and $d$, 
 %  $ \vert CN(\sigma^w, d)\vert={w\choose d}.$
 % \end{corollary}
 % \noindent {\bf Proof.}
 % Substituting $s$ with 2 in Eqn. (\ref{eqn4}),  we have
 % $$\vert CN(\sigma^w, d)\vert ={w-d-1 \choose 0} + 
 % {w-d\choose 1} 
 % + \cdots + {w-1\choose d} = {w-d-1 +d+1\choose d}={w\choose d}.$$
 % $\Box$
%\todo[inline]{Should we change w to n, as w denotes words earlier? }

Using a similar argument, we can prove the following formula for the size of the super-condensed $d$-neighborhood for unary words. 
\begin{proposition}
  %  Let $\Sigma$ be an alphabet consisting of $s$ characters and  $\sigma \in \Sigma$. 
  Let $\Sigma$ be an alphabet consisting of $s$ characters, 
  { and let $\sigma \in \Sigma$.
    For any { integers $n$ and $d$ such that} $0 < d { <  n-1}$,}
    \begin{equation}
    \vert SCN({ \sigma^n}, d)\vert  = \sum_{n-d\leq m\leq n}
        {m-2\choose d+m-n}(s-1)^{d+m-n}.
    \end{equation}
    % \begin{eqnarray}
    % \label{eqn4_super}
    %  &&   \vert SCN(w, d)\vert  \nonumber \\ 
    %    & =&  \sum_{n-d\leq m\leq n}
    %     {m-2\choose d+m-n}(s-1)^{d+m-n}.
    % \end{eqnarray}
    In particular, if $s=2$,
    \begin{equation}
        \label{eqn4_super_binary}
        \vert SCN({ \sigma^n}, d)\vert = {n-1\choose d}.
    \end{equation}
\end{proposition}
%\todo[inline]{the bound n is replaced with n-2 to answer the question of Reviewer3}
%% ----------------------------------------------------------------------------------------------
\section{Condensed Neighborhood for Arbitrary Words}
\label{sec:arbitrary}

Our main result in this section is  a novel upper bound  on the sizes of the condensed neighborhoods.  This bound leads to a proof of a conjecture introduced in~\cite{DBLP:conf/stringology/ChauveMP21}. By definition, it is also an upper bound on the size of the super-condensed neighborhood for arbitrary words.
%\todo[inline]{change W to w, its length to n.}
 
\begin{proposition}  
\label{prop3}
Let  $w=w_1w_2\cdots w_n$ be a word over an alphabet $\Sigma$ and $n>0$, i.e., $w\in \Sigma^{+}$,  and  { let $d$ be an integer such that $0{ < } d< n$. }
For any  $x \in CN(w, d)$,  we have that (i) $d_{lev}(x, w)=d$,  and (ii) { in any optimal alignment between $x$ and $w$, $x_m$ belongs to a match column.
}
\end{proposition}
\noindent {\bf Proof}. 
Let $x=x_1x_2\cdots x_m\in CN(w, d)$.
 { Since the Levenshtein distance satisfies the triangle inequality, } if $d_{lev}(x, w)\leq d-1$ then
\begin{eqnarray*}
  &&  d_{lev}(x_1x_2\cdots x_{m-1}, w) \\
  &\leq & d_{lev}(x_1x_2\cdots x_{m-1}, x)+ d_{lev}(x, w)\\
  &\leq & 1 + (d-1)\leq d.
\end{eqnarray*}
This implies that $x_1x_2\cdots x_{m-1}\in N(w, d)$, contradicting $x\in CN(w, d)$. 

Let $A$ be an optimal alignment of $w$ and $x$ that consists of $t$ columns where $x$ appears in the first row and $w$  in the second row. 
Assume that $x_m$ appears in a mismatch or deletion column ${a_{i} \brack b_i}$ (so $a_i=x_m$).
%where $i\geq m$ and $b_i\neq x_m$ in $\Sigma$ or $b_i=-$. 
Then,  $A$ has the following structure:
\begin{center}
\begin{tabular}{cccccccc}
  $a_1$& $a_2$ & $\cdots$ & $a_{i-1}$
  & $x_m$ & $-$ & $\cdots$ & $-$\\
  $b_1$ & $b_2$ & $\cdots$ & $b_{i-1}$ &
  $b_i$ & $b_{i+1}$ & $\cdots$ & $b_{t}$
\end{tabular}
\end{center}
where $i\geq m$ and $b_i = -$ or $b_i\in\Sigma \setminus\{x_m\}$.
{ Then, $x'=x_1x_2\cdots x_{m-1}\in N(w, d)$, contradicting $x\in CN(w, d)$.
Indeed, if $b_i=-$, removing from $A$ the column ${x_m\brack b_i}$ results in an alignment between $x'$ and $w$ of cost $d_{lev}(x, w)-1$, while, if $b_i\neq -$, replacing $a_i$ by $-$ in $A$ results in an alignment between $x'$ and $w$ of cost $d_{lev}(x, w)$.}
% the following alignment
% between $x_1x_2\cdots x_{m-1}$ and $w$:
% \begin{center}
% \begin{tabular}{cccccccc}
%   $a_1$& $a_2$ & $\cdots$ & $a_{i-1}$
%   & $-$ & $-$ & $\cdots$ & $-$\\
%   $b_1$ & $b_2$ & $\cdots$ & $b_{i-1}$ &
%   $b_i$ & $b_{i+1}$ & $\cdots$ & $b_{t}$
% \end{tabular}
% \end{center}
% has score $d_{lev}(x, w)$ (if $b_i \neq -$) or $d_{lev}(x, w)-1$ (if $b_i=-$).
% This implies that $x_1x_2\cdots x_{m-1}\in N(w, d)$, contradicting $x\in CN(w, d)$.
$\qed$
\\

\begin{definition}
Let ${A}$ be an optimal alignment between two words $x$ and $y$ with $k$ match/mismatch columns with indices, { ordered increasingly}, $(i_1, i_2, \cdots, i_k)$. 
${  A}$ 
{ is said }
to be the leftmost optimal alignment { between $x$ and $y$} if, for any other optimal alignment ${  B}$ between $x$ and $y$, ${  B}$ has $\ell\geq k$ match/mismatch columns and { the { increasing sequence of the indices $(p_1, p_2, \cdots, p_\ell)$ of these columns} of ${  B}$,  is lexicographically greater than $(i_1, i_2, \cdots, i_k)$, that is, there exists $t$ such that $i_j=p_j$ for each $j\leq t$ and $i_{t+1}<p_{t+1}$. 
}
\end{definition}
%\\

%\subsection{New Counting Formulas}

\begin{lemma}
\label{lemma1}
Let $w=w_1w_2\cdots w_n$ and $x=x_1x_2\cdots x_m$ be two words over { an alphabet} $\Sigma$ such that $x\in CN(w, d)$,
{ where $0<d<n$.  }
Let ${  A}$ be the leftmost optimal alignment between $x$ and $w$,  where $x$ appears in the first row.
If $w_j$  is the last character { of $w$ not belonging to} an insertion column 
%${- \brack w_j}$ 
in ${A}$, then 
(i)  the { last character $x_m$ of $x$ and $w_j$ form a match column ${x_m\brack w_j}$},   
and 
(ii) { the column 
%${x_{m-1} \brack -}$ { right} 
in $A$ immediately before ${x_m\brack w_j}$ is not a deletion column.}
\end{lemma}
\noindent {\bf Proof.} 
By Proposition~\ref{prop3}, $x_m$ appears in a match { column in $A$.} 
As $w_j$ is the last character of $w$ that is not in an insertion column, { $w_j$ forms a match column with $x_m$, which proves (i)}. 

{ To prove (ii), assume that} ${x_{m-1} \brack -}$ appears { immediately } before ${x_m\brack w_j}$ in $A$.
Then replacing these two columns by the column ${x_{m-1} \brack w_j}$ results in an alignment between $x_1x_2\cdots x_{m-1}$ and $w$ that also contains $d$ insertion, deletion and mismatch columns. 
This implies that $x_1x_2\cdots x_{m-1}\in N(W, d)$, { contradicting that $x\in CN(w, d)$.}
$\qed$
\\
%\todo[inline]{For consistency, in Proof of  lemma 1,  we assume x always appears in the first row in an alignment between x and w.In Prop 4, W is changed to w, the len is denoted by n}

\begin{proposition}
    \label{prop:upper}
    Let $w$ be a word of length $n$ over an alphabet $\Sigma$ such that $\vert\Sigma \vert=s$. 
    Then, for any $d$ such that $0< d < n$, 
    {\small 
    \begin{eqnarray}
    \label{Ineq6}
        \vert CN(w, d) \vert  & \leq & \sum_{0\leq i\leq d} {n \choose i} (s-1)^{d-i}  \nonumber  \times \sum_{0\leq j\leq d-i} { n-i-1\choose j} {n+d-2i-2j-{ 2} \choose d-i-j}.
    \end{eqnarray}
    }
\end{proposition}
\noindent {\bf Proof.} 
%Let $w$ be a word on an alphabet $\Sigma$ such that $|w|=n$ and $|\Sigma|=s$. 
For each word $x\in CN(w, d)$, we consider the leftmost optimal alignment ${A}$ between $x$ and $w$, in which $x$ appears in the first row and $w$ appears in the second row. 
{ Assume that} $A$ contains:
\begin{itemize}
\item $i$ insertion columns for some $i$ such that $0\leq i\leq d$, 
\item $j$ mismatch columns for some $j$ such that $0\leq j\leq d-i$, 
\item $d-i-j$ deletion columns,
\item $n-i-j\geq 1$ match columns.
\end{itemize}
There are ${n \choose i}$  possible ways of selecting the $i$ { characters of $w$ that belong to insertion columns}. { Once these inserted positions are fixed, by Lemma~\ref{lemma1}.(i), the last character of $w$ that is not in any of these positions must form a match column with the last character of $x$.
So this match column can be followed in $A$ only by insertion columns.}
Therefore,  there are ${n-i-1\choose j}$ possible ways of selecting $j$ { positions in $w$ that belong to mismatch columns, and the $n-i-j$ remaining positions of $w$ belong to match columns.
To complete $A$ we need to decide where are blocks of consecutive deletions (deletion blocks), which defines the structure of $A$, and which characters to assign to $x$ in mismatch and deletion columns.

In any optimal alignment, a deletion column cannot occur immediately before an insertion column as otherwise both columns could be combined to form a match or a mismatch column in an alignment of lower cost than the cost of $A$.
In $A$, a deletion block can not appear immediately before a mismatch column, as otherwise one could shift the character of $w$ in this mismatch column to the column immediately to its left, to obtain an alignment of cost at most the cost of $A$, and this would contradict that $A$ is the leftmost optimal alignment between $x$ and $w$. 
By Lemma~\ref{lemma1}.(ii), a deletion block can not occur immediately before the last match column.
So the $d-i-j$ deletion columns are divided into $n-i-j-1$ blocks (some of which could be empty) each occurring immediately before one of the $n-i-j-1$ first match columns. 
There are $n-i-j-1+d-i-j-1\choose d-i-j$ possible ways to split $d-i-j$ deletions into $n-i-j-1$ blocks, and each such configuration defines a unique alignment structure that needs to by completed by assigning a character to each of the positions of $x$ that participate to a deletion or a mismatch column.

Positions assigned to mismatch columns can be assigned $s-1$ possible characters.
By construction, a deletion block is always followed by a match column ${x_i \brack w_j}$ (say $x_i=w_j=\sigma$); if the character $x_i$ appears in the deletion block, then the match column could be replaced by a deletion, with $w_j$ being shifted to the left to align with any occurrence of $\sigma$ in the deletion block, which contradicts that $A$ is the leftmost optimal alignment.
So the $j$ mismatch columns and the $d-i-j$ deletion columns are restricted to $s-1$ possible characters each.

Combining these facts, we obtain the right-hand side of Inequality~(\ref{Ineq6}).

Last, the alignment structure defined above might not be the structure of a leftmost optimal alignment.
For example a deletion block could follow immediately an insertion column, in which case both columns could be combined into a match or mismatch column to define an alignment of lower cost. 
Another configuration that is not compatible with a leftmost optimal alignment would be the case where $w_i=\sigma$ belongs to an insertion column followed by a match column with $w_{i+1}=\sigma$.
This is why we do not have equality between the right-hand side of Inequality~(\ref{Ineq6}) and $\vert CN(w,d)\vert$.}
$\qed$
%% ----------------------------------------------------------------------------------------------
\section{A Simple Upper Bound Formula for the Condensed Neighborhood}
\label{sec:conjecture}

We finally prove an elegant upper bound for the size of the condensed neighborhood for arbitrary words 
{ 
%(Inequality~(\ref{Ineq49})), 
that was}
conjectured in~\cite{DBLP:conf/stringology/ChauveMP21}, whose values { can however be} a few times larger than the { upper-bound given in}  Proposition~\ref{prop:upper} { as illustrated, for $s=2$ and small values of $n$ and $d$, in} Table~\ref{table1}\footnote{The code used to generate Table~\ref{table1} is available at \url{https://github.com/cchauve/CondensedNeighbourhoods/tree/ARXIV2025}}.

\begin{theorem}
    \label{thm:conjecture}
    For any  word $w$ of length $n$ over an alphabet $\Sigma$ such that $\vert\Sigma\vert =s$, and any  { $0< d < n$}
    \begin{eqnarray}
    \label{Ineq49}
        \vert {CN}(w, d)\vert \leq \frac{(2s-1)^d n^d}{d!}.
    \end{eqnarray}
\end{theorem}

\begin{table*}[htbp]
\begin{center}
   \caption{{ The values of (Top) upper-bound~(\ref{Ineq6}) and (Bottom) upper-bound~(\ref{Ineq49}) for $s=2$, $n=4, 6, 8, 10$ and $1\leq d\leq n-1$.}}
    \label{table1}
        \begin{tabular}{c|rrrrrrrrr}
            \hline
            \textbf{$|w| \setminus d$} &\textbf{1} & \textbf{2} & \textbf{3} & \textbf{4} & \textbf{5} & \textbf{6} & \textbf{7} &  \textbf{8} & \textbf{9}\\
            \hline
            \hline
            \textbf{4}   & 10 & 37  & 63    &        &        &         &         &         & \\
            \textbf{6}   & 16 & 108 & 403   & 935    & 1,526  &         &         &         & \\	
            \textbf{8}	 & 22 & 215 & 1,235 & 4,678  & 12,587 & 25,943  & 44,936  &         & \\
            \textbf{ 10} & 28 & 358 & 2,775 & 14,638 & 56,168 & 164,969 & 389,994 & 784,085 & 1,414,039\\
            \hline
            \hline
            \textbf{4}   & 12 & 72  & 288   &        &          &           &           &            & \\
            \textbf{6}   & 18 & 162	& 972	& 4,374	 & 15,746	&           &           &            & \\	
            \textbf{8}	 & 24 & 288	& 2,304	& 13,824 & 66,355   & 265,420	& 910,014   &            & \\
            \textbf{ 10} & 30 & 450	& 4,500	& 33,750 & 202,500	& 1,012,500	& 4,339,285	& 16,272,321 & 54,241,071 \\
        \hline
        \end{tabular}
\end{center}
\end{table*}

\begin{lemma}
    \label{lemma5}
    { For integers $n$, $d$ and $j$ such that $0\leq j \leq d\leq n$,} 
    \begin{eqnarray}
        {n\choose j} { n +d -2j \choose d-j} &\leq & \frac{(n+d/2+1/2)^{d-j}(n-d/2+3/2)^j}{j!(d-j)!}.
    \end{eqnarray}
\end{lemma} 
\noindent {\bf Proof.} Provided in Appendix.

\medskip\noindent {\bf Proof of Theorem~\ref{thm:conjecture}.}
{ 
First, for any $n \geq 1$, and $0\leq d<n$, by expanding $(2s-1)^d$ as $(1+2(s-1))^d$ with the Binomial Theorem, we have
\begin{eqnarray}
    \label{eq1.1}
    \frac{(2s-1)^d n^d}{d!} &= & \sum_{0\leq x\leq d} (s-1)^{d-x}  \frac{2^{d-x}n^d}{x!(d-x)!}.
\end{eqnarray}
Next,  for integers $n$, $d$ and $x$ such that $0\leq x\leq d<n$, 
\begin{eqnarray}
    {n \choose x}  \sum_{0\leq j\leq d-x} {n-x-1\choose j} { n +d -x-2j-1 \choose d-x-j}&\leq &  \frac{2^{d-x}n^d}{x!(d-x)!}.
    \label{eq1.2}
\end{eqnarray}
To prove Inequality~(\ref{eq1.2}), replacing $d-x$ with $b$ and $n-1$ with $m$ (note that $b \geq 0$ and $m \geq 0$), we have
\begin{eqnarray*} 
    {n \choose x}  \sum_{0\leq j\leq b} {n-x-1\choose j} { n +b-2j-1 \choose b-j} & \leq & \frac{(m+1)^x}{x!} 
    \left(\prod_{k=0}^{x-1}\frac{m+1-k}{m+1} \right) 
    \sum_{0\leq j \leq b} 
    \left[
    \left(\prod_{k=0}^{j-1} \frac{m-x-k}{m-k}\right){m \choose j} {m+b-2j \choose b-j}\right]\\
    &\leq &  \frac{(m+1)^x}{x!} \sum_{0\leq j\leq b} {m \choose j} {m+b-2j \choose b-j}.
\end{eqnarray*}
Therefore, Inequality~(\ref{eq1.2}) holds if, %for $b\leq n-1$,
\begin{eqnarray}
    \label{eq1.3}
    \sum_{0\leq j\leq b} {m\choose j} {m +b -2j \choose b-j} 
    & \leq & \frac{2^b(m+1)^b}{b!}
\end{eqnarray}
Inequality~(\ref{eq1.3}) follows from Lemma~\ref{lemma5} and the Binomial Theorem expansion
\[
\frac{2^b(m+1)^b}{b!} = \sum_{0\leq j\leq b}\frac{(m+b/2+1/2)^{b-j}(m-b/2+3/2)^j}{j!(b-j)!}.
\]
Last, by Inequality~(\ref{eq1.1}) and Inequality~(\ref{eq1.2})
\begin{eqnarray*}
    \frac{(2s-1)^d n^d}{d!}&\geq &  \sum_{0\leq x\leq d} 
    (s-1)^{d-x}{n\choose x} \times \sum_{0\leq j\leq d-x}{n-x-1\choose j}{n+d-x-2j-1\choose d-x-j}
\end{eqnarray*}
Proposition~\ref{prop:upper}, together with
\[
{n+d-2x-2j-1\choose d-x-j} \geq {n+d-2x-2j-2\choose d-x-j}
\]
proves the Theorem.}
$\qed$

%    &\geq &  \sum_{0\leq x\leq d} 
%     (s-1)^{d-x}{n\choose x}  \\
%     && \times \sum_{0\leq j\leq d-x}{n-x-1\choose j}{n+d-2x-2j-{ 2}\choose d-x-j} \\
%     && \;\;\; \mbox{{ (Decrease the top part of the last binomial coeff.)}}\\ 
%     &\geq& \vert CN(w, d)\vert. \;\;\;\;\;\;
%     \mbox{(Proposition~\ref{prop:upper})} 
%     \qed
% \end{eqnarray*}
%% ----------------------------------------------------------------------------------------------
\section{Conclusion}
\label{sec:conclusion}

{ 
In this note, we proved several enumerative results on word neighborhoods: exact formulas for the size of the condensed and super-condensed neighborhoods of unary words, and two upper-bounds for the size of the condensed neighborhood of arbitrary words. 
These results suggest several avenues for further research.}

Our results on the size of condensed and super-condensed neighborhoods of unary words extend the result introduced in~\cite{DBLP:conf/cpm/Charalampopoulos20b} for whole neighborhoods. 
It was also shown in~\cite{DBLP:conf/cpm/Charalampopoulos20b} that unary words have the smallest neighborhoods among all words of the same length over a given alphabet, thus leading to lower bounds for the size of neighborhoods.
It is thus natural to ask if a similar property holds for condensed and super condensed neighborhoods, namely that unary words have the smallest condensed or super condensed neighborhoods. 

\begin{figure}
    \label{fig:fig1}
    \begin{center}
        \includegraphics[width=18cm]{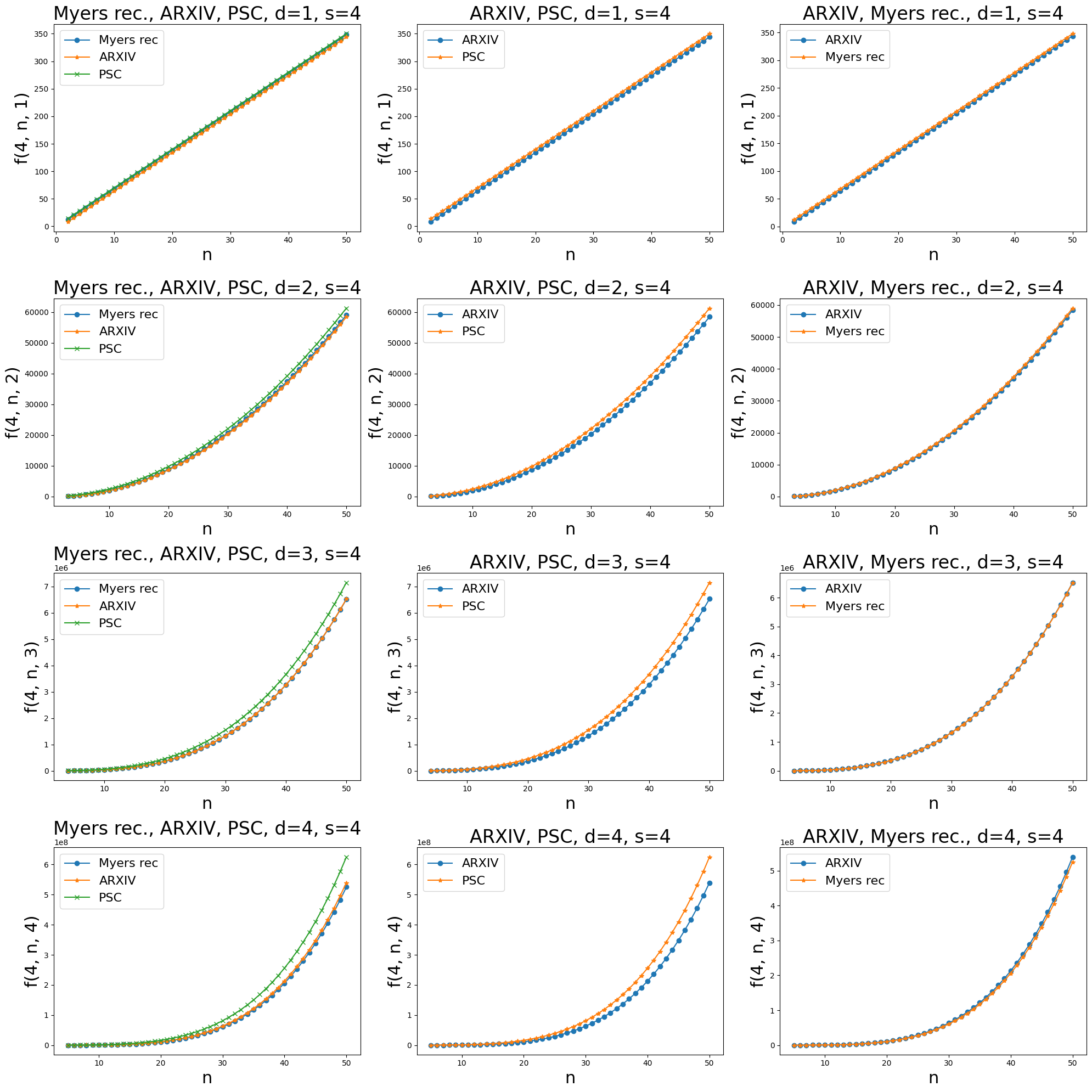}
    \end{center}
    \caption{Upper-bounds on the size of the condensed neighborhood for an alphabet of size $s=4$, words of length up to $n=50$ and Levensthein distance $d=1,2,3,4$. \textbf{Myers rec.}: upper-bound defined by the recurrences described in~\cite{DBLP:books/daglib/p/Myers13}. \textbf{ARXIV}: Proposition~\ref{prop:upper}. \textbf{PSC}: Theorem~\ref{thm:conjecture}.
    The code to generate this figure is available at
    \url{https://github.com/cchauve/CondensedNeighbourhoods/tree/ARXIV2025}.}
\end{figure}

{ 
The novel upper-bound on the size of the condensed neighborhood for arbitrary words that we introduce in Proposition~\ref{prop:upper} is based on counting alignments of words in the condensed neighborhood, similar to what was done in~\cite{DBLP:books/daglib/p/Myers13}.
However, while in~\cite{DBLP:books/daglib/p/Myers13} such alignments were counted through a set of recurrences, we provide a non-recursive formula based on a double summation. 
In both cases, the proposed formulas are upper-bounds as several counted alignments can define the same word of the condensed neighborhood.
Experiments (Appendix Fig.~\ref{fig:fig1}) show that both upper-bounds are very close with no pattern allowing to claim that one is always better than the other one.
It remains open to design improved recurrences or formulas for counting alignments of words in the condensed neighborhood that excludes more alignments defining the same word. 
}

The conjecture we proved in Theorem~\ref{thm:conjecture} was introduced in~\cite{DBLP:conf/stringology/ChauveMP21}, where its simple form allowed to use it to suggest, through experimental results, that it could lead to a slightly larger window { of average-case linear time complexity for the approximate pattern matching used in BLAST~\cite{DBLP:books/daglib/p/Myers13} compared to the analysis based on the recurrences introduced in~\cite{DBLP:books/daglib/p/Myers13}.}
It remains open to use Theorem~\ref{thm:conjecture} in a theoretically rigorous average-case complexity analysis of the BLAST algorithm.
The comparison between the two upper bounds we proved (Table~\ref{table1} { and Appendix Fig.~\ref{fig:fig1}}) also shows that the more complex upper bound introduced in Proposition~\ref{prop:upper} is much sharper than the upper bound of Theorem~\ref{thm:conjecture}.
{ Given its relatively simple form, it is open to investigate whether it is amenable to be used to analyze the average-case time complexity of the BLAST algorithm.}

%% ----------------------------------------------------------------------------------------------
\bibliographystyle{elsarticle-num}
\bibliography{references}

\vfill\newpage
\section*{Appendix}

\begin{lemma}
    \label{lemma3a}
    { For integers $n$ and $d$ such that $0<d<n$, and any $t$ such that $d > t\geq \frac{d}{4}$},
    \begin{eqnarray*}
        (n-t+1)^2\leq 
        \left(n+\frac{d}{2}+\frac{1}{2}\right)\left(n - \frac{d}{2}+\frac{3}{2}\right).
    \label{eqn777}
    \end{eqnarray*}
\end{lemma}
\noindent {\bf Proof.}
% If $t=\frac{d}{4}$, then $1\leq d\leq n$ and thus
% \begin{eqnarray*}
%     && \left(n+\frac{d}{2}+\frac{1}{2}\right)\left(n-\frac{d}{2}+\frac{3}{2}\right) - \left(n-\frac{d}{4}+ 1\right)^2\\
%     & =&  -\frac{1}{4} + d - \frac{5}{16} d^2 + \frac{1}{2}d n \geq - \frac{5}{16} d^2 + \frac{1}{2}d n \geq 0.
% \end{eqnarray*}
% For any $n\geq d$ and $ d> t> \frac{d}{4} \geq 0$,
% \begin{eqnarray*}
%     (n-t+1)^2 < \left(n -\frac{d}{4}+1\right)^2,
% \end{eqnarray*}
% and thus Inequality~(\ref{eqn777}) holds.
{ As $\frac{d}{4}\leq t < d$, 
$(n-t+1)^2 \leq (n-d/4+1)^2$. Next,
\begin{eqnarray*}
    && \left(n+\frac{d}{2}+\frac{1}{2}\right)\left(n-\frac{d}{2}+\frac{3}{2}\right) - \left(n-\frac{d}{4}+ 1\right)^2\\
    & =&  -\frac{1}{4} + d - \frac{5}{16} d^2 + \frac{1}{2}d n.
\end{eqnarray*}
As $d\geq 1$ and $n>d$, the expression above is always positive which proves the Lemma.
$\qed$}

\begin{lemma}
    \label{lemma4}
    { For integers $n$ and $d$ such that $0<d<n$, and any $t$ such that $0\leq t\leq \frac{d}{4}$,}
    \begin{eqnarray*}
        &&  (n + 2t)(n + 2t-1)(n -t + 1)^2 
        \nonumber \\ 
        &\leq & \left(n+\frac{d}{2}+\frac{1}{2}\right)^3\left(n-\frac{d}{2}+\frac{3}{2}\right). 
        %\label{eqn_31}
    \end{eqnarray*}    
\end{lemma}
\noindent {\bf Proof.} 
For { $0\leq t\leq d/4$, }
\begin{eqnarray*}
    && (n+d/2+1/2)^3(n-d/2+3/2) 
    \\
    && \;\;- (n+2t)(n+2t-1)(n-t+1)^2  \nonumber  \\
    &=& 2n^3+ 4n^2+\frac{9}{4}n + \frac{3}{16}-\frac{1}{16}d^4 - \frac{1}{4}d^3 n + \frac{3}{4} d^2 n \nonumber \\
    && + \frac{3}{8} d^2 + d n^3 + 3 d n^2 + \frac{9}{4} d n+\frac{1}{2} d \nonumber  \\
    && - 4 t^4 + 4 t^3 n + 10 t^3 + 3 t^2 n^2 - 3 t^2 n - 8 t^2 \\
    && - 2 t n^3 - 6 t n^2 - 2 t n + 2 t  \nonumber  \\  
    &\geq &     2n^3+ 4n^2+\frac{9}{4}n + \frac{3}{16}-\frac{1}{16}d^4 - \frac{1}{4}d^3 n + \frac{3}{4} d^2 n \\
    && + \frac{3}{8} d^2 + d n^3 + 3 d n^2 + \frac{9}{4} d n \nonumber  \\
    && - 4 t^4  - 3 t^2 n - 8 t^2 - 2 t n^3 - 6 t n^2 - 2 t n  \\
    && \;\;\;\;\;(\mbox{Delete $d/2$ and { all} positive terms containing $t$})  \nonumber\\  
    &\geq & 2n^3+ 4n^2+\frac{9}{4}n + \frac{3}{16}-\frac{1}{16}d^4 - \frac{1}{4}d^3 n + \frac{3}{4} d^2 n \\
    && + \frac{3}{8} d^2 + d n^3 + 3 d n^2 + \frac{9}{4} d n \nonumber\\
    && - \frac{1}{64} d^4  - \frac{3}{16} d^2 n - \frac{1}{2} d^2 - \frac{1}{2} d n^3 - { \frac{6}{4}  dn^2} - \frac{1}{2}d n  \\
    && \;\;\;\;\;(\mbox{Substitute  $t$ with $d/4$})     \nonumber   \\
    &=& 2n^3+ { 4n^2} +\frac{9}{4}n + \frac{3}{16}-\frac{5}{64}d^4 - \frac{1}{4}d^3 n + \frac{9}{16} d^2 n \\
    && -\frac{1}{8} d^2 +\frac{1}{2} d n^3 + { \frac{6}{4}  dn^2} + \frac{7}{4} d n  \nonumber  \\
    &=& 2n^3+ \left(4n^2-\frac{1}{8}d^2\right)+\frac{9}{4}n + \frac{3}{16} + \frac{9}{16} d^2 n \\
    && +\left( d n^3  -\frac{5}{64}d^4 - \frac{1}{4}d^3n\right) + { \frac{6}{4}  dn^2} + \frac{7}{4} d n   \nonumber  \\
    &\geq & 0. \;\;\;\;\; (\mbox{From $n\geq d$}) \nonumber \qed
\end{eqnarray*}  
% Inequality~(\ref{eqn_31}) is proved.
% $\qed$

\medskip
\noindent {\bf Proof of Lemma~\ref{lemma5}.} 
If $j=0$, { it is straightforward to verify that the inequality holds. So we consider now that $j>0$.}
% \begin{eqnarray*}
%     &&  {n\choose j} { n +d -2j \choose d-j} \\
%     &=&   \frac{ (n + d)(n+d-1)\cdots (n+1)}{j!(d-j)!} \\
%     & \leq &  \frac{ (n + d/2+1/2)^{d-0}}{0!(d-0)!},
% \end{eqnarray*}

\noindent ({\bf Case 1})~~
For any integer $j$ such that $0<j <  d/2$, define $k=d/2-j$. (Note that $k$ is not an integer if $d$ is odd, { but $2k$ is always an integer}.) Then, we have:
$$d/2=j+k, \;\; d-j=d/2+k=j+2k,$$
and
\begin{eqnarray*}
    X  &=&   j!(d-j)! {n\choose j} { n +d -2j \choose d-j} \\
%    &=& n(n-1)\cdots (n-j+1)\\
%    && \times (n+2k)(n+2k-1)\cdots (n-j+1)\\
    &=&(n+2k)\cdots (n+1) \left(n(n-1)\cdots (n-j+1)\right)^2.
\end{eqnarray*}
%We consider the following two cases: 
\noindent ({\bf Case 1.1})~~ Assume $j\leq d/4\leq k$. 
Then, we have,
\begin{eqnarray*}
    && (n+2k)(n+2k-1)\cdots (n+1)\\
    &=& { \left(\prod_{i=1}^{2(k-j)}(n+2j+1)\right)(n+2j)\cdots (n+1)}\\
%    &=& (n+2k-1)\cdots (n+2j+1)\\
%    && \times (n+2j)\cdots (n+1)\\
    &\leq & (n+j+k+1/2)^{2(k-j)}(n+2j)\cdots (n+1)\\
    &=& \left(n+\frac{d}{2}+\frac{1}{2}\right)^{2k-2j} (n+2j)\cdots (n+1), 
\end{eqnarray*}
{ We rewrite $(n+2j)\cdots (n+1)\left(n(n-1)\cdots (n-j+1)\right)^2$ as
\[
\prod_{i=0}^{j-1}(n+2(j-i))(n+2(j-i)-1)(n-(j-i)+1)^2
\]}
and by Lemma~\ref{lemma4}, { which applies as $0\leq j-i\leq d/4$,} 
\begin{eqnarray*}
%    X &\leq & \left(n+\frac{d}{2}+\frac{1}{2}\right)^{2k-2j} (n+2j)\cdots (n+1) \\
%    && \times \left(n(n-1)\cdots (n-j+1)\right)^2\\
    X &\leq & \left(n+\frac{d}{2}+\frac{1}{2}\right)^{2k-2j} \left(n+\frac{d}{2}+\frac{1}{2}\right)^{3j}
    \times \left(n-\frac{d}{2}+\frac{3}{2}\right)^{j}\\
    &=& \left(n+\frac{d}{2}+\frac{1}{2}\right)^{d-j}\left(n -\frac{d}{2}+\frac{3}{2}\right)^{j}.
\end{eqnarray*}
The inequality is proved for $j\leq d/4\leq k.$

\noindent ({\bf Case 1.2})~~ Assume  $k\leq d/4 \leq j$. By Lemma~\ref{lemma4}, we have
\begin{eqnarray*}
    && (n+2 \lfloor k\rfloor)\cdots (n+1)(n(n-1)\cdots (n-\lfloor k\rfloor +1))^2\\
    & \leq&  \left(n+\frac{d}{2}+\frac{1}{2}\right)^{3\lfloor k\rfloor }\left(n-\frac{d}{2}+\frac{3}{2}\right)^{\lfloor k\rfloor},
\end{eqnarray*}
and, if $k$ is not an integer,  
$ 
n+2k \leq n+d/2 \leq n+d/2+1/2.
$
Therefore,
\begin{eqnarray}
    && (n+2 k)\cdots(n+1)(n(n-1)\cdots (n-\lfloor k\rfloor +1))^2 \nonumber \\
    &\leq &  \left(n+\frac{d}{2}+\frac{1}{2}\right)^{3 \lfloor k\rfloor +1 }\left(n-\frac{d}{2}+\frac{3}{2}\right)^{\lfloor k\rfloor}. 
    \label{eqn8}
\end{eqnarray}
If $\lfloor k\rfloor + j$ is even, as $j+k=d/2$ and $j\geq d/4$, by Lemma~\ref{lemma3a},
\begin{eqnarray}
    && (n-\lfloor k\rfloor )^2(n-\lfloor k\rfloor -1)^2...(n-j+1)^2  \nonumber \\
    &=& \left[\prod^{(\lfloor k\rfloor +j)/2-1}_{t=\lfloor k\rfloor} (n-t)(n-j+t+1)\right]^2  \nonumber \\
    &\leq & \left[\prod^{(\lfloor k\rfloor +j)/2-1}_{t=\lfloor k\rfloor} (n-j+1/2)^2\right]^2  \nonumber \\
    &=& \left(n+\frac{d}{2}+\frac{1}{2}\right)^{j-\lfloor k\rfloor}\left(n-\frac{d}{2}+\frac{3}{2}\right)^{j-\lfloor k\rfloor}.
    \label{eqn9}
\end{eqnarray}

If $\lfloor k\rfloor +j$ is odd,  
by Lemma~\ref{lemma3a}, 
%
%\[(w-j+1)^2 \leq (w+d/2+1/2)(w-d/2+3/2)\]
%and thus,
\begin{eqnarray}
    && (n-\lfloor k\rfloor )^2(n-\lfloor k\rfloor -1)^2\cdots (n-(j-1)+1)^2  \nonumber  \times (n-j+1)^2 \nonumber \\
    &\leq & \left(n+\frac{d}{2}+\frac{1}{2}\right)^{j-1-\lfloor k\rfloor}\left(n-\frac{d}{2}+\frac{3}{2}\right)^{j-1-\lfloor k\rfloor} \nonumber  \times (n-j+1)^2 \nonumber \\
    &\leq & \left(n+\frac{d}{2}+\frac{1}{2}\right)^{j-\lfloor k\rfloor}\left(n-\frac{d}{2}+\frac{3}{2}\right)^{j-\lfloor k\rfloor}.
    \label{eqn10}
\end{eqnarray}

By Inequalities~(\ref{eqn8})-(\ref{eqn10}), 
we have proved that 
\begin{eqnarray*} 
    X &\leq &  \left(n+\frac{d}{2}+\frac{1}{2}\right)^{j+2k}\left(n-\frac{d}{2}+\frac{3}{2}\right)^j \\
    &= & \left(n+\frac{d}{2}+\frac{1}{2}\right)^{d-j}\left(n-\frac{d}{2}+\frac{3}{2}\right)^j,
\end{eqnarray*}
as $d=2j+2k$. 

\noindent {\bf (Case 2}) Let $j=d-k>d/2$ for some $k<d/2$.
\begin{eqnarray*}
    X &=&    j!(d-j)! {n\choose d-k} { n -d+2k \choose k} \\
    &=&    n(n-1)\cdots (n-k+1) (n-k)\cdots (n-d/2)\\
    && \times (n-d/2-1)\cdots (n-d+k+1) \\
    && \times (n-d+2k) (n-d+2k-1)\cdots (n-d+k+1)\\
    & \leq & n(n-1)\cdots(n-k+1) \left(n-\frac{d}{2}+\frac{3}{2}\right)^{d-2k}\\
    && \times (n-d+2k)(n-d+2k-1)\cdots (n-d+k+1)\\
    &=& \left(n-\frac{d}{2}+\frac{3}{2}\right)^{d-2k}\;\;\prod^{k-1}_{i=0}[(n-i)(n-d+k+i+1)] \\
    &\leq & \left(n-\frac{d}{2}+\frac{3}{2}\right)^{d-2k}\;\;\prod^{k-1}_{i=0}\left(n-\frac{d}{2}+\frac{k}{2}+\frac{1}{2}\right)^2\\
    &\leq &  \left(n+\frac{d}{2}+\frac{1}{2}\right)^k\left(n-\frac{d}{2}+\frac{3}{2}\right)^{k+d-2k}\\
    &=& \left(n+\frac{d}{2}+\frac{1}{2}\right)^{d-j}\left(n-\frac{d}{2}+\frac{3}{2}\right)^{j},
\end{eqnarray*}
where the last 
inequality follows from Lemma~\ref{lemma3a} and the fact that
$\frac{d-k}{2}=\frac{j}{2} >\frac{d}{4}.$
$\qed$

\vfill~\newpage

\onecolumn

\end{document}